\documentstyle[12pt]{article}
\advance\hoffset by -7mm  
\makeatletter
\@addtoreset{equation}{section}
\makeatother

\renewcommand{\title}[1]{\null\vspace{25mm}

\noindent{\Large{\bf #1}}\vspace{10mm}

\noindent {\large By }}
\newcommand{\authors}[1]{\noindent{\large #1}\vspace{3mm}

}
\newcommand{\address}[1]{\noindent #1\vspace{5mm}

}
\renewcommand{\abstract}[1]{\vspace{19mm}

\noindent{\small{\em Abstract.} #1}\vspace{2mm}

} 
\font\tenof=msym10 at 12pt
\def\R{\mbox{\tenof R}}
\def\Z{\mbox{\tenof Z}}
\newtheorem{proposition}{Proposition}
\newcommand{\ap}{a^{\dagger}}
\newcommand{\tlambda}{\tilde{\lambda}}
\newcommand{\tmu}{\tilde{\mu}}
\newcommand{\tnu}{\tilde{\nu}}

\begin{document}
%
%
\title{Representation Theory of Deformed Oscillator Algebras}
\authors{Christiane Quesne \footnote{Directeur de recherches FNRS; E-mail:
cquesne@ulb.ac.be} and Nicolas Vansteenkiste \footnote{E-mail: nvsteen@ulb.ac.be}}
\address{Physique Nucl\'eaire Th\'eorique et Physique Math\'ematique, 
Universit\'e Libre de Bruxelles, Campus de la Plaine CP229, Boulevard du 
Triomphe, B-1050 Brussels, Belgium}
%
%
\abstract{ 
The representation theory of deformed oscillator algebras, defined in terms of an
arbitrary function of the number operator~$N$, is developed in terms of the
eigenvalues of a Casimir operator~$C$. It is shown that according to the nature of
the $N$ spectrum, their unitary irreducible representations may fall into one out of
four classes, some of which contain bosonic, fermionic or parafermionic
Fock-space representations as special cases. The general theory is illustrated by
classifying the unitary irreducible representations of the Arik-Coon,
Chaturvedi-Srinivasan, and Tamm-Dancoff oscillator algebras, which may be
derived from the boson one by the recursive minimal-deformation procedure of
Katriel and Quesne. The effects on non-Fock-space representations of the minimal
deformation and of the quommutator-commutator transformation, considered in
such a procedure, are studied in detail.}
%
%
\section{Introduction}
Since the pioneering works of Arik and Coon~\cite{arik},
Kuryshkin~\cite{kuryshkin}, Biedenharn~\cite{biedenharn}, and
Macfarlane~\cite{macfarlane}, many forms of deformed oscillator algebras have
been considered and played an important role in the construction of
$q$-deformed Lie algebras (see e.g.~\cite{hayashi,fairlie}). They have found various
applications to physical problems, such as the description of systems with
non-standard statistics~\cite{greenberg,cqa}, the construction of integrable
lattice models~\cite{bogoliubov}, the algebraic treatment of quantum mechanical
exactly solvable systems~\cite{daska92,cqb}, of pairing correlations in nuclear
physics~\cite{bonatsos}, and of vibrational spectra of diatomic and polyatomic
molecules~\cite{chang}, as well as the search for nonlinearities due to
high-intensity electromagnetic fields~\cite{manko}.\par
%
%
The necessity to introduce some order in the rich and varied choice of deformed
commutation relations did however appear soon and various classification
schemes were therefore proposed~\cite{jannussis,daska91,solomon,meljanac}.
More recently, a unifying recursive procedure was introduced, generating at
appropriate steps all the familiar deformed oscillators, along with some
multiparametric generalizations~\cite{katriel}. Each iteration consists in a
minimal deformation of a commutator into a quommutator, followed by a
transformation of the latter into a new commutator, to which it is equivalent
within the corresponding Fock space.\par
%
%
As it was already observed in Ref.~\cite{katriel}, the equivalence between
quommutators and corresponding commutators, which is a central ingredient of the
recursive minimal deformation procedure, is not valid any more in the additional
non-Fock-space representations, which are known to exist in general for deformed
oscillator algebras. Although various works have been devoted to determining such
representations for some particular
algebras~\cite{kulish,rideau,chaichian,kosinski}, there still remains a need for
a general theory.\par
%
%
The purpose of the present paper is twofold: first to fill in this gap by discussing
the representation theory of general deformed oscillator algebras; then to
illustrate both the effects of minimal deformation and of the
quommutator-commutator transformation on the non-Fock-space representations
by studying some selected examples.\par
%
%
The recursive minimal deformation procedure is briefly reviewed in Sec.~2. The
representation theory of general deformed oscillator algebras is then developed in
Sec.~3, and illustrated on some examples in Sec.~4. Finally, Sec.~5 contains the
conclusion.\par
%
%
\section{Recursively Minimally-Deformed Oscillators}
Let us consider a given oscillator algebra~${\cal A}_0$, generated by the operators
$N = N^{\dagger}$, $\ap$, $a = \left(\ap\right)^{\dagger}$, satisfying the
commutation relations
\begin{eqnarray}
  \left[N, \ap\right] & = & \ap, \qquad \left[N, a \right] = - a,    \label{eq:2.1} \\
  \left[a , \ap\right] & = & f_0(N),     \label{eq:2.2}
\end{eqnarray}
for some function $f_0(N) = \left(f_0(N)\right)^{\dagger}$. The algebra ${\cal A}_0$
will serve as a starting point for a recursive procedure, wherein other oscillator
algebras will be generated~\cite{katriel}.\par
%
%
Let us assume that at the $k$th step, we have obtained an algebra~${\cal A}_k$,
still generated by $N$, $\ap$,~$a$, and satisfying Eq.~(\ref{eq:2.1}), but with
Eq.~(\ref{eq:2.2}) replaced by
\begin{equation}
  \left[a, \ap\right] = f_k(N),      \label{eq:2.3}
\end{equation}
where $f_k(N) = \left(f_k(N)\right)^{\dagger}$. Then the next minimal
deformation~$\tilde{\cal A}_k$ of this algebra is defined by Eq.~(\ref{eq:2.1})
and
\begin{equation}
  \left[a, \ap\right]_{q_{k+1}} = f_k(N),     \label{eq:2.4}
\end{equation}
where $q_{k+1}$ is some real parameter, and the left-hand side of
Eq.~(\ref{eq:2.4}) is a quommutator, defined by $\left[a, \ap\right]_{q_{k+1}} \equiv
a \ap - q_{k+1} \ap a$.\par
%
%
The minimally-deformed relation~(\ref{eq:2.4}) implies that in the bosonic
Fock-space representation, i.e., with respect to the eigenvectors $|n\rangle$ of the
number operator~$N$, corresponding to the eigenvalues $n = 0$, 1, 2,~$\ldots$, the
operators $\ap$ and~$a$ satisfy the relations
\begin{equation}
  \ap |n\rangle = \sqrt{F_{k+1}(n+1)}\, |n+1\rangle, \qquad
  a |n\rangle = \sqrt{F_{k+1}(n)}\, |n-1\rangle,     \label{eq:2.5}
\end{equation}
where the vacuum state~$|0\rangle$ is assumed to fulfil the condition
\begin{equation}
  a |0\rangle = 0,      \label{eq:2.5a}
\end{equation}
and the function $F_{k+1}(n)$ is defined by
\begin{equation}
  F_{k+1}(n) = \sum_{i=0}^{n-1} q_{k+1}^i f_k(n-1-i).    \label{eq:2.6}
\end{equation}
In Eq.~(\ref{eq:2.6}), $\sum_{i=0}^{-1} \equiv 0$ so that $F_{k+1}(0) = 0$ in
accordance with Eqs.~(\ref{eq:2.5}) and~(\ref{eq:2.5a}). The corresponding function
of the number operator $F_{k+1}(N)$, which satisfies the equation
\begin{equation}
  F_{k+1}(N+1) - q_{k+1} F_{k+1}(N) = f_k(N),    \label{eq:2.6a}
\end{equation}
is referred to as the structure function of the algebra~$\tilde{\cal A}_k$.\par
%
%
It follows that the algebra~${\cal A}_{k+1}$, defined by Eq.~(\ref{eq:2.1}) and
\begin{equation}
  \left[a, \ap\right] = f_{k+1}(N),      \label{eq:2.7}
\end{equation}
where
\begin{equation}
  f_{k+1}(N) = \left(f_{k+1}(N)\right)^{\dagger} \equiv F_{k+1}(N+1) - F_{k+1}(N),
  \label{eq:2.8}
\end{equation}
is equivalent to~$\tilde{\cal A}_k$ in such a Fock-space representation. In other
words, both algebras~$\tilde{\cal A}_k$ and~${\cal A}_{k+1}$ have the same
structure function~$F_{k+1}(N)$.\par
%
%
By iterating the transformations ${\cal A}_k \to \tilde{\cal A}_k \to {\cal
A}_{k+1}$, one gets a sequence of deformed oscillator algebras depending upon an
increasing number of parameters.\par
%
%
Both ${\cal A}_k$ and $\tilde{\cal A}_k$ have a nonvanishing central element or
Casimir operator, defined\linebreak by~\cite{katriel,nvsteen}
\begin{equation}
  C_k = F_k(N) - \ap a,     \label{eq:2.9}
\end{equation}
and
\begin{equation}
  \tilde C_k = q_{k+1}^{-N} \left(F_{k+1}(N) - \ap a\right) = q_{k+1}^{-N} C_{k+1},
  \label{eq:2.10}
\end{equation}
respectively, which may be used to characterize their irreducible representations.
As in the Fock-space representation, $F(N) |0\rangle = a |0\rangle = 0$, it follows
from (\ref{eq:2.9}) and (\ref{eq:2.10}) that $C_k$ and $\tilde C_k$ have a
vanishing eigenvalue in such a representation.\par
%
%
In Sec.~4, we shall consider as examples the first few iterations obtained by
starting from the standard boson oscillator algebra, for which
\begin{equation}
  f_0(N) = 1, \qquad F_0(N) = N, \qquad C_0 = N - \ap a.    \label{eq:2.10a}
\end{equation}
In such a case, the first minimal deformation~$\tilde{\cal A}_0$ is the Arik-Coon
oscillator algebra~\cite{arik,kuryshkin}, for which
\begin{equation}
  \left[a, \ap\right]_{q_1} = 1.     \label{eq:2.11}
\end{equation}
Its structure function and Casimir operator are given by
\begin{equation}
  F_1(N) = [N]_{q_1} \equiv \frac{q_1^N - 1}{q_1 -1}, \qquad \tilde C_0 = 
  \frac{1 - q_1^{-N}}{q_1 - 1} - q_1^{-N} \ap a,      \label{eq:2.12}
\end{equation}
respectively. The Fock-space equivalent algebra~${\cal A}_1$ corresponds to the
Chaturvedi-Srinivasan oscillator~\cite{chaturvedi}, for which
\begin{equation}
  \left[a, \ap\right] = q_1^N, \qquad F_1(N) = [N]_{q_1}, \qquad C_1 = [N]_{q_1}
  - \ap a.        \label{eq:2.13}
\end{equation}
The second minimal deformation~$\tilde{\cal A}_1$ is the Chakrabarti-Jagannathan
two-parameter oscillator algebra~\cite{chakrabarti}, for which
$\left[a, \ap\right]_{q_2} = q_1^{N}$. It has two important special cases: the
Biedenharn~\cite{biedenharn} and Macfarlane~\cite{macfarlane} oscillator algebra,
corresponding to $q_2 = q_1^{-1}$, and the Tamm-Dancoff oscillator
algebra~\cite{odaka}, for which $q_2 = q_1$, and
\begin{equation}
  \left[a, \ap\right]_{q_1} = q_1^{N}, \qquad F_2(N) = q_1^{N-1} N, \qquad
  \tilde C_1 = q_1^{-1} N - q_1^{-N} \ap a.        \label{eq:2.14}
\end{equation}
\par
%
%
Many more examples can be found in Ref.~\cite{katriel}.\par
%
%
\section{Representation Theory of Deformed Oscillator Algebras}
In the present section, we shall review some general properties of the unitary
irreducible representations (unirreps) of the deformed oscillator algebras
considered in the previous one. For such purpose, it is enough to consider the case
of~$\tilde{\cal A}_k$, as that of~${\cal A}_k$ can be obtained from it by
restricting the $q_{k+1}$ values to $q_{k+1}=1$. For simplicity's sake, we shall
define the algebra commutation relations by Eq.~(\ref{eq:2.1}) and
\begin{equation}
  \left[a, \ap\right]_q = f(N) = \left(f(N)\right)^{\dagger},        \label{eq:3.1}
\end{equation}
and denote the corresponding structure function and Casimir operator by~$F(N)$
and~$C$, respectively. From (\ref{eq:2.6}) and~(\ref{eq:2.10}), the latter satisfy
the relations
\begin{equation}
  F(N+1) - q F(N) = f(N), \qquad C = q^{-N} \left(F(N) - \ap a\right),     \label{eq:3.2}
\end{equation}
from which it follows that
\begin{equation}
  \ap a = F(N) - q^N C, \qquad a \ap = F(N+1) - q^{N+1} C.              \label{eq:3.3}
\end{equation}
We shall classify the unirreps of this algebra under the assumption that the
spectrum of~$N$ is discrete and nondegenerate.\par
%
%
Let us start with a normalized simultaneous eigenvector $|c, \nu_0\rangle$ of the
Casimir operator~$C$, defined in~(\ref{eq:3.2}), and of~$N$, corresponding to the
eigenvalues $c$ and~$\nu_0$ respectively,
\begin{equation}
  C |c, \nu_0\rangle = c |c, \nu_0\rangle, \qquad N |c, \nu_0\rangle = \nu_0
  |c, \nu_0\rangle, \qquad \langle c, \nu_0 |c, \nu_0\rangle = 1.    \label{eq:3.4}
\end{equation}
From (\ref{eq:2.1}), it results that the vectors
\begin{equation}
  |c, \nu_0 + n) = \left\{\begin{array}{ll}
                                     \left(\ap\right)^n |c, \nu_0\rangle & \mbox{if
                                          $n=1,2,\ldots$,} \\[0.2cm]
                                     a^{-n} |c, \nu_0\rangle & \mbox{if $n=-1,-2,\ldots$,}  
                                   \end{array}\right.    \label{eq:3.5}
\end{equation}
are also simultaneous eigenvectors of $C$ and~$N$,
\begin{equation}
  C |c, \nu_0 + n) = c |c, \nu_0 + n), \qquad N |c, \nu_0 + n) = (\nu_0 + n) 
  |c, \nu_0 + n),      \label{eq:3.6}
\end{equation}
as long as they are nonvanishing. In~(\ref{eq:3.5}), we use a round bracket instead of
an angular one to denote unnormalized states. Definition~(\ref{eq:3.5}) can be
extended to $n=0$ by setting $|c, \nu_0) = |c, \nu_0\rangle$.\par
%
%
Eq.~(\ref{eq:3.3}) implies that
\begin{equation}
  \ap a |c, \nu_0 + n) = \lambda_n |c, \nu_0 + n), \qquad a \ap |c, \nu_0 + n) = \mu_n
  |c, \nu_0 + n),         \label{eq:3.7}
\end{equation}
where
\begin{equation}
  \lambda_n = F(\nu_0+n) - q^{\nu_0+n} c, \qquad \mu_n  = F(\nu_0+n+1) 
  - q^{\nu_0+n+1} c = \lambda_{n+1}.                 \label{eq:3.8}
\end{equation}
As eigenvalues of a positive operator, only those $\lambda_n$ that are positive or
null are admissible in a unitary representation. In particular, the condition
\begin{equation}
  \lambda_0 = F(\nu_0) - q^{\nu_0} c \ge 0          \label{eq:3.9}
\end{equation}
restricts the possible values of $c$ and~$\nu_0$.\par
%
%
So it is straightforward to derive the following unitarity conditions:
\begin{proposition}
If there exists some $m_1 \in \{-1, -2, -3, \ldots\}$ such that $\lambda_{m_1} < 0$,
and $\lambda_n \ge 0$ for $n = 0$, $-1$, $\ldots$,~$m_1+1$, then an irreducible
representation of a deformed oscillator algebra can be unitary only if
$\lambda_{n_1} = 0$ for some $n_1 \in \{0, -1, \ldots, m_1+1\}$. If there exists
some $m_2 \in \{2, 3, 4, \ldots,\}$ such that $\lambda_{m_2} < 0$, and $\lambda_n
\ge 0$ for $n = 0$, $1$, $\ldots$,~$m_2-1$, then it can be unitary only if
$\lambda_{n_2} = 0$ for some $n_2 \in \{1, 2, \ldots, m_2-1\}$.
\end{proposition}
\noindent{\it Proof.} In the first part of the proposition, we must have $|c,
\nu_0+m_1) = 0$ as otherwise $\ap a$ would have a negative eigenvalue. This
implies that $a |c, \nu_0+m_1+1) = 0$, which can be achieved in two ways, either 
$|c, \nu_0+m_1+1) = 0$, or $|c, \nu_0+m_1+1) \ne 0$ and $\lambda_{m_1+1} = 0$.
In the former case, we can proceed in the same way and find that at least one of
the conditions $\lambda_{m_1+2} = 0$, $\lambda_{m_1+3} = 0$, $\ldots$,
$\lambda_{-1} = 0$, or $|c, \nu_0-1) = 0$ must be satisfied. But the last one is
equivalent to $\lambda_0 = 0$, since $|c, \nu_0) \ne 0$ by hypothesis. This
concludes the proof of the first part of the proposition. The second part can be
demonstrated in a similar way by using $a \ap$, and $\mu_n = \lambda_{n+1}$
instead of $\ap a$, and~$\lambda_n$.\par
%
%
According to Proposition~1, the unirreps may belong to one out of four classes. If
there exists some $n_1 \in \{0, -1, -2, \ldots\}$ such that $\lambda_{n_1} = 0$, and
$\lambda_n > 0$ for $n = n_1+1$, $n_1+2$, $\ldots$, 0, 1, 2,~$\ldots$, then $|c,
\nu_0+n_1)$ satisfies the relation $a |c, \nu_0+n_1) = 0$, and is an eigenvector of
$\ap a$ and~$N$ with eigenvalues $\lambda_{n_1} = 0$ and $\tnu_0 = \nu_0 + n_1$.
By repeating construction~(\ref{eq:3.5}) with $|c, \nu_0\rangle$, $\nu_0$,
$\lambda_0$ replaced by $|c, \tnu_0\rangle$, $\tnu_0$, $\tlambda_0 = 0$,
respectively , we obtain for the corresponding normalized states~\footnote{In
Eq.~(\ref{eq:3.10}), we assume that $\prod_{i=1}^0 \equiv 1$. A similar convention
is used in subsequent formulae too.}
\begin{equation}
  |c, \tnu_0 + n\rangle = \left(\prod_{i=1}^n \tlambda_i\right)^{-1/2} 
  \left(\ap\right)^n |c, \tnu_0\rangle, \qquad n = 0, 1, 2, \ldots,    \label{eq:3.10}
\end{equation}
where $\tlambda_n = \lambda_{n+n_1} = F(\tnu_0+n) - q^n F(\tnu_0)$, and the
Casimir operator eigenvalue~$c$ is entirely determined by~$\tnu_0$ through the
relation $c = q^{-\tnu_0} F(\tnu_0)$. The states~(\ref{eq:3.10}) carry an
infinite-dimensional unirrep, characterized by a lower bound~$\tnu_0$ (bounded
from below or BFB~unirrep). In basis~(\ref{eq:3.10}), the generators are represented
by
\begin{eqnarray}
  a |c, \tnu_0 + n\rangle & = & \sqrt{\tlambda_n}\, |c, \tnu_0 + n - 1\rangle, \qquad
            \ap |c, \tnu_0 + n\rangle = \sqrt{\tlambda_{n+1}}\, |c, \tnu_0 + n + 1
            \rangle, \nonumber \\[0.2cm]
  N |c, \tnu_0 + n\rangle & = & \left(\tnu_0 + n\right) |c, \tnu_0 + n\rangle.
            \label{eq:3.11}     
\end{eqnarray}
In the special case where $\tnu_0 = 0$, we obtain a bosonic Fock-space
representation of type~(\ref{eq:2.5}), wherein the spectrum of~$N$ is $\{0, 1, 2,
\ldots\}$, and $c = F(0) = 0$.\par
%
%
If, on the contrary, there exists some $n_2 \in \{1, 2, 3, \ldots\}$ such that
$\lambda_{n_2} = 0$, and $\lambda_n > 0$ for $n = n_2-1$, $n_2-2$,~$\ldots$, 0,
$-1$, $-2$,~$\ldots$, then $|c, \nu_0+n_2-1)$ satisfies the relation $\ap |c,
\nu_0+n_2-1) = 0$, and is an eigenvector of $a\ap$ and~$N$ with eigenvalues
$\mu_{n_2-1} = \lambda_{n_2} = 0$ and $\tnu_0 = \nu_0 + n_2 - 1$. If we repeat
construction~(\ref{eq:3.5}) by starting from $|c, \tnu_0\rangle$, $\tnu_0$,
$\tmu_0 = \tlambda_1 = 0$, instead of $|c, \nu_0\rangle$, $\nu_0$, $\mu_0 =
\lambda_1$, we obtain for the corresponding normalized states
\begin{equation}
  |c, \tnu_0 + n\rangle = \left(\prod_{i=0}^{|n|-1} \tlambda_{-i}\right)^{-1/2} 
  a^{-n} |c, \tnu_0\rangle, \qquad n = 0, -1, -2, \ldots,    \label{eq:3.12}
\end{equation}
where $\tlambda_n = \lambda_{n+n_2-1} = F(\tnu_0+n) - q^{n-1} F(\tnu_0+1)$, and
$c$ is again determined by~$\tnu_0$ through the relation $c = q^{-\tnu_0-1}
F(\tnu_0+1)$. Such states now carry an infinite-dimensional unirrep, characterized
by an upper bound~$\tnu_0$ (bounded from above or BFA unirrep). The
representation of the generators in basis~(\ref{eq:3.12}) is still given
by~(\ref{eq:3.11}), but where $n$ now takes the values indicated in~(\ref{eq:3.12}),
instead of those shown in~(\ref{eq:3.10}).\par
%
%
It may also happen that there exist both $n_1 \in \{0, -1, -2, \ldots\}$, and $n_2
\in \{1, 2, 3, \ldots\}$ such that $\lambda_{n_1} = \lambda_{n_2} = 0$, and
$\lambda_n > 0$ for $n = n_1+1$, $n_1+2$,~$\ldots$, $-1$, 0, 1,~$\ldots$, $n_2-2$,
$n_2-1$. The corresponding unirrep is then finite-dimensional (FD~unirrep), and
may be characterized by its lower and upper bounds, $\tnu_0 = \nu_0 + n_1$ and
$\tnu_0 + n_2 - n_1 - 1 = \nu_0 + n_2 - 1$, or alternatively by $\tnu_0$ and $p =
n_2 - n_1 - 1$. It is spanned by the $d = p + 1$ normalized states
\begin{equation}
  |c, \tnu_0 + n\rangle = \left(\prod_{i=1}^n \tlambda_i\right)^{-1/2} 
  \left(\ap\right)^n |c, \tnu_0\rangle, \qquad n = 0, 1, \ldots, p,    \label{eq:3.13}
\end{equation}
where $\tlambda_n = \lambda_{n+n_1} = F(\tnu_0+n) - q^n F(\tnu_0)$, and $c =
q^{-\tnu_0} F(\tnu_0) = q^{-\tnu_0-p-1} F(\tnu_0+p+1)$. They still satisfy
Eq.~(\ref{eq:3.11}), but we note that now
\begin{equation}
  \ap |c, \tnu_0 + p\rangle = 0.       \label{eq:3.14}
\end{equation}
The unirrep is an order-$p$ parafermionic Fock-space representation if $c = \tnu_0
= 0$. It is fermionic in the special case where $p=1$.\par
%
%
Finally, if $\lambda_n > 0$ for $n \in \Z$, we get an unbounded unirrep
(UB~unirrep), which may be characterized by~$c$ and by the fractional
part~$\tnu_0$ of~$\nu_0$ (i.e., $\nu_0 = [\nu_0] + \tnu_0$, where $0 \le \tnu_0 <
1$, and $[\nu_0]$ denotes the largest integer contained in~$\nu_0$), as different
values of~$[\nu_0]$ lead to equivalent unirreps. Its representation space is
spanned by the states
\begin{eqnarray}
  |c, \tnu_0+n\rangle &= & \left(\prod_{i=1}^n \tlambda_i\right)^{-1/2}
            \left(\ap\right)^n |c, \tnu_0\rangle, \qquad n=0,1,2,\ldots, \nonumber \\
  |c, \tnu_0+n\rangle & = & \left(\prod_{i=0}^{|n|-1} \tlambda_{-i}\right)^{-1/2}
            a^{-n} |c, \tnu_0\rangle, \qquad n=-1,-2,\ldots,         \label{eq:3.15} 
\end{eqnarray} 
where $\tlambda_n = \lambda_{n-[\nu_0]}$, and the generators are still
represented by Eq.~(\ref{eq:3.11}).\par
%
%
\section{Some Selected Examples}
In the present section, we shall apply the theory developed in the previous one to
some of the deformed oscillator algebras considered in Sec.~2. The results are
summarized in Tables~1, 2, and~3.\par
%
%
As explained in Sec.~2, the starting algebra ${\cal A}_0$ of the recursive procedure
considered here is the boson oscillator algebra, defined by Eqs.~(\ref{eq:2.1}),
(\ref{eq:2.2}), and~(\ref{eq:2.10a}). It is worth emphasizing that contrary to the
Heisenberg algebra for which the number operator is defined as $N \equiv \ap a$,
the boson oscillator algebra has some non-Fock-space representations.
From~(\ref{eq:3.8}) and~(\ref{eq:3.9}), we indeed obtain that
\begin{equation}
  \lambda_n = \nu_0 + n - c, \qquad \mbox{where } \nu_0 \ge c,    \label{eq:4.1}
\end{equation}
may become negative for $n < c - \nu_0$. Hence unitarity imposes that there exists
some $n_1 \in \{0, -1, -2, \ldots\}$ such that $\lambda_{n_1} = \nu_0 + n_1 - c =
0$. The algebra has therefore BFB~unirreps, characterized by
\begin{equation}
  \tnu_0 = \nu_0 + n_1 = c, \qquad \tlambda_n = n,     \label{eq:4.2}
\end{equation}
where $\tnu_0$ may take any real value. For $\tnu_0 = c \ne 0$, such
representations are non-Fock-space unirreps.\par
%
%
We shall successively review the cases where $0 < q \ne 1$, and $q < 0$. The latter
is omitted in most studies, because the corresponding algebras are considered as
deformations of the fermion oscillator algebra, instead of the boson one. It is
worth noting however that in some definitions of deformed oscillator
algebras~\cite{daska91}, both a commutation and an anticommutation relations are
assumed. We chose here to keep only one of them. As explained in
Refs.~\cite{nvsteen,oh}, such a modified definition leads to the existence of a
Casimir operator. Considering negative $q$ values for the minimally-deformed
oscillator algebras is therefore in some way equivalent to selecting
anticommutation relations instead of the commutation ones associated with
positive $q$ values.\par
%
%
Note that for $q<0$, except when otherwise stated, we shall restrict $\nu_0$ to
integer values so that $q^{\nu_0}$ is well defined.\par
%
%
\subsection{The Arik-Coon Oscillator Algebra}
\subsubsection{Positive Values of the Deforming Parameter}
For the Arik-Coon oscillator algebra $\tilde{\cal A}_0$~\cite{arik,kuryshkin},
defined by Eqs.~(\ref{eq:2.1}) and~(\ref{eq:2.11}), we find from (\ref{eq:2.12}),
(\ref{eq:3.8}), and~(\ref{eq:3.9}) that for~$q>0$
\begin{equation}
  \lambda_n = \left(\frac{1}{q-1} - c\right) q^{\nu_0+n} - \frac{1}{q-1}, \qquad
  \mbox{where } c\le \frac{1 - q^{-\nu_0}}{q-1},      \label{eq:4.3}
\end{equation}
may be an increasing, constant or decreasing function of~$n$ according to the
values taken by~$q$ and~$c$. To classify its unirreps, we have to distinguish
between the cases where $0<q<1$ and $q>1$.\par
%
%
Whenever $0<q<1$, we note from Eq.~(\ref{eq:4.3}) that the Casimir operator
eigenvalue may satisfy either of the conditions $c\le (q-1)^{-1}$, or
$(q-1)^{-1} < c \le \left(1 - q^{-\nu_0}\right)/(q-1)$. In the former case, $\lambda_n
> 0$ for any $n \in \Z$, so we obtain UB~unirreps, whereas in the latter case,
$\lambda_n$ may become negative for some negative $n$~values, so we get
BFB~unirreps, characterized by
\begin{equation}
  \tnu_0 = \nu_0 + n_1, \qquad c = \frac{1 - q^{-\tnu_0}}{q-1}, \qquad \tlambda_n
  = [n]_q.              \label{eq:4.4}
\end{equation}
Here $[n]_q$ is defined as in Eq.~(\ref{eq:2.12}), and $n_1 \in \{0, -1, -2, \ldots\}$,
so that $\tnu_0$ may take any real value. Note that for $c = (q-1)^{-1}$, the
UB~unirrep degenerates into a unirrep for which
\begin{equation}
  \ap a = a \ap = (1-q)^{-1}.       \label{eq:4.5}
\end{equation}
\par
%
%
Whenever $q>1$, we always have $c \le \left(1 - q^{-\nu_0}\right)/(q-1) <
(q-1)^{-1}$. Since $\lambda_n$ may again become negative for some negative
$n$~values, we obtain BFB~unirreps, similar to those defined in~(\ref{eq:4.4}).\par
%
%
In the limit where $q \to 1^-$ or~$1^+$, the only surviving unirreps are the
BFB~ones, which go over into those of~${\cal A}_0$, defined in~(\ref{eq:4.2}). The
UB~unirreps, which diverge for $q \to 1^-$, are referred to as classically singular
representations~\cite{aizawa}.\par
%
%
Our results do agree with those previously derived by Kulish~\cite{kulish} by a
similar type of approach. The method used here, as well as in Ref.~\cite{kulish},
contrasts with that of Chaichian {\it et al.}~\cite{chaichian}. Indeed the latter do
not postulate the existence of a number operator, hence of Eq.~(\ref{eq:2.1}). Their
unirrep classification is therefore not performed in terms of a Casimir
operator~$C$, but in terms of some noncentral element $K \equiv a\ap - \ap a$,
whose sign cannot change in a given unirrep. Whenever $K\ne 0$, they set $|K| =
q^M$, where the operators $M$, $\ap$, and~$a$ satisfy some relations similar to
Eq.~(\ref{eq:2.1}). The connection between their approach~\footnote{It is worth
noting that Chaichian {\it et al.} call any BFB~unirrep a Fock-space representation,
whereas we do reserve this name for a very specific BFB~unirrep.} and ours is
easily established by noting that $K$~can be rewritten in terms of our
operators~$N$ and~$C$ as $K = q^N (1 + (1-q) C)$. Hence $K>0$, $K=0$, and $K<0$
correspond to $c > (q-1)^{-1}$ if $0<q<1$, or $c < (q-1)^{-1}$ if $q>1$, $c =
(q-1)^{-1}$, and $c < (q-1)^{-1}$ if $0<q<1$, respectively, and for~$K\ne0$, one may
set $M = N + \log_q|1 + (1-q) c|$.\par
%
%
\begin{table}[htb]

\caption{Unirrep classification for the Arik-Coon oscillator algebra. The cases
where $q=1$ and $q=-1$ correspond to the boson and fermion oscillator algebras,
respectively.}

\vspace{1cm}
\begin{tabular}{lll}          
  \hline\\[-0.2cm] 
  $q$ & Type & Characterization \rule[-0.3cm]{0cm}{0.6cm}\\[0.2cm]
  \hline\\[-0.2cm] 
  $q>1$       & BFB & $\tnu_0\in\R$, $c=q^{-\tnu_0} [\tnu_0]_q$, $\tlambda_n = [n]_q$
                       \\[0.2cm]
  $q=1$      & BFB & $\tnu_0\in\R$, $c=\tnu_0$, $\tlambda_n = n$ \\[0.2cm]
  $0<q<1$   & BFB & $\tnu_0\in\R$, $c=q^{-\tnu_0} [\tnu_0]_q$, $\tlambda_n =
                       [n]_q$\\[0.2cm]
                 & UB   & $0\le\tnu_0<1$, $c\le (q-1)^{-1}$, $\tlambda_n = [\tnu_0+n]_q -
                       c q^{\tnu_0+n}$\\[0.2cm]
  $-1<q<0$ & BFB & $\tnu_0\in\Z$, $c=q^{-\tnu_0} [\tnu_0]_q$, $\tlambda_n =
                       [n]_q$\\[0.2cm]
                 & UB  & $0\le\tnu_0<1$, $c=(q-1)^{-1}$, $\tlambda_n=(1-q)^{-1}$ 
                       \\[0.2cm]
  $q=-1$    & FD  & $\tnu_0\in 2\Z$, $p=1$, $c=0$, $\tlambda_n = \left(1 - 
                       (-1)^n\right)/2$ \\[0.2cm]
                & FD   & $\tnu_0\in 2\Z + 1$, $p=1$, $c=-1$, $\tlambda_n = \left(1 - 
                       (-1)^n\right)/2$ \\[0.2cm]
                & UB   & $\tnu_0=0$, $-1<c<-1/2$ or $-1/2<c<0$, \\[0.2cm] 
                &        & $\tlambda_n = (-1)^{n+1}c + \left(1 - (-1)^n\right)/2$ \\[0.2cm]
                & UB   & $0\le\tnu_0<1$, $c=-1/2$, $\tlambda_n=1/2$  \\[0.2cm]
  $q<-1$   & BFA  & $\tnu_0\in\Z$, $c=q^{-\tnu_0-1} [\tnu_0+1]_q$, $\tlambda_n =
                       [n-1]_q$  \\[0.2cm]
               & UB   & $0\le\tnu_0<1$, $c=(q-1)^{-1}$, $\tlambda_n = (1-q)^{-1}$
\\[0.2cm]
  \hline 
\end{tabular}
\end{table}
%
%
\subsubsection{Negative Values of the Deforming Parameter}
From Eq.~(\ref{eq:4.3}), it is obvious that for any negative $q$~value, and $c =
(q-1)^{-1}$, there exists a degenerate UB~unirrep, for which Eq.~(\ref{eq:4.5}) is
valid, and which may be characterized by any $\tnu_0$ such that $0\le \tnu_0
<1$.\par
%
%
Assuming now $c\ne (q-1)^{-1}$ and $\nu_0 \in \Z$, Eq.~(\ref{eq:4.3}) becomes
\begin{equation}
  \lambda_n = (-1)^{\nu_0+n+1} \left(\frac{1}{1+|q|} + c\right) |q|^{\nu_0+n}
  + \frac{1}{1+|q|},               \label{eq:4.6}
\end{equation}
where
\begin{eqnarray}
  c & \le & \frac{|q|^{-\nu_0}-1}{1+|q|} \qquad \mbox{if } \nu_0 \in 2\Z, 
               \nonumber \\
  c & \ge & - \frac{|q|^{-\nu_0}+1}{1+|q|} \qquad \mbox{if } \nu_0 \in 2\Z+1.
               \label{eq:4.7}
\end{eqnarray}
For successive $n$~values, $\lambda_n$ oscillates around the positive constant $(1
+ |q|)^{-1}$. To classify the unirreps, we have to distinguish between the cases
where $0<|q|<1$, $|q|>1$, and $|q|=1$.\par
%
%
Whenever $0<|q|<1$, $|\lambda_n|$ decreases from~$+\infty$ to~$(1+|q|)^{-1}$.
Hence if $(-1)^{\nu_0} ((1+|q|)^{-1} + c) > 0$, $\lambda_n$ may become
negative for some negative even $n$~values. Unitarity then imposes that there
exists some $n_1 \in \{0, -2, -4, \ldots\}$ such that
\begin{equation}
  \lambda_{n_1-2} < 0, \qquad \lambda_{n_1} = 0, \qquad \lambda_{n_1-1},
  \lambda_{n_1+1}, \lambda_{n_1+2}, \ldots >0,      \label{eq:4.9}
\end{equation}
corresponding to $m_1 = n_1-2$ in Proposition~1. We therefore obtain
BFB~unirreps, characterized by
\begin{equation}
  \tnu_0 = \nu_0 + n_1 \in 2\Z, \qquad c = \frac{|q|^{-\tnu_0}-1}{1+|q|}, \qquad
           \tlambda_n = \frac{1+(-1)^{n+1}|q|^n}{1+|q|},            \label{eq:4.10}
\end{equation}
if $\nu_0 \in 2\Z$, or
\begin{equation}
  \tnu_0 = \nu_0 + n_1 \in 2\Z + 1, \qquad c = - \frac{|q|^{-\tnu_0}+1}{1+|q|},
           \qquad \tlambda_n = \frac{1+(-1)^{n+1}|q|^n}{1+|q|},    \label{eq:4.11} 
\end{equation}
if $\nu_0 \in 2\Z + 1$. If, on the contrary, $(-1)^{\nu_0} ((1+|q|)^{-1} + c) < 0$,
$\lambda_n$ may become negative for some negative odd $n$~values. This shows
that Eq.~(\ref{eq:4.9}) must be satisfied for some $n_1 \in \{-1, -3, -5, \ldots\}$.
Hence, we again obtain BFB~unirreps, but this time the $\tnu_0$, $c$, and
$\tlambda_n$~values that characterize them are given by Eq.~(\ref{eq:4.10}) or
(\ref{eq:4.11}) according to whether $\nu_0 \in 2\Z+1$ or $\nu_0 \in 2\Z$. The
results obtained in the various cases can be put together by writing that for
$0<|q|<1$, there exist BFB~unirreps specified by
\begin{equation}
  \tnu_0 \in \Z, \qquad c = \frac{q^{-\tnu_0}-1}{1-q}, \qquad \tlambda_n = [n]_q,
  \label{eq:4.12}
\end{equation}
where the $\tnu_0$~value is arbitrary.\par
%
%
Whenever $|q|>1$, $|\lambda_n|$ increases from $(1+|q|)^{-1}$ to~$+\infty$, so that
$\lambda_n$ may become negative for some positive even or odd values, according
to whether $(-1)^{\nu_0} ((1+|q|)^{-1} + c)$ is positive or negative. Unitarity now
imposes that there exists some $n_2 \in \{2, 4, 6, \ldots\}$ or  $n_2 \in \{1, 3, 5,
\ldots\}$ respectively, such that
\begin{equation}
  \lambda_{n_2+2} < 0, \qquad \lambda_{n_2} = 0, \qquad \lambda_{n_2+1},
  \lambda_{n_2-1}, \lambda_{n_2-2}, \ldots >0,      \label{eq:4.13}
\end{equation}
corresponding to $m_2 = n_2+2$ in Proposition~1. By proceeding as in the case
where $0<|q|<1$, we conclude that there exist BFA~unirreps, characterized by
\begin{equation}
  \tnu_0 \in \Z, \qquad c = \frac{q^{-\tnu_0-1}-1}{1-q}, \qquad \tlambda_n =
  [n-1]_q,        \label{eq:4.14}
\end{equation}
for any $\tnu_0$ value.\par
%
%
Finally, for $|q|=1$, corresponding to the fermion oscillator algebra,
Eq.~(\ref{eq:4.6}) becomes
\begin{equation}
  \lambda_n = (-1)^{\nu_0+n+1} \left(c + \frac{1}{2}\right) + \frac{1}{2},
  \label{eq:4.15}        
\end{equation}
where $c\le 0$ or $c\ge -1$ according to whether $\nu_0 \in 2\Z$ or $\nu_0 \in
2\Z+1$. In the former case, we obtain that for $c=0$, $\lambda_n = 0$ for any $n
\in 2\Z$, while $\lambda_n = 1$ for any $n \in 2\Z+1$, thereby showing that there
exist FD~unirreps, characterized by $p=1$ and any $\tnu_0 \in 2\Z$. For $\tnu_0 =
0$, this is the standard fermionic Fock-space representation. For $c=-1$,
$\lambda_n = 1$ for any $n \in 2\Z$, while $\lambda_n = 0$ for any $n \in 2\Z+1$.
So we again get FD~unirreps characterized by $p=1$, but this time $\tnu_0 \in
2\Z+1$. They can be derived from the previous ones by interchanging the roles of
$\ap$ and~$a$. For $-1<c<-1/2$ or $-1/2<c<0$, $\lambda_n$ is always positive,
hence the corresponding unirreps are UB~ones. Similar results are obtained by
starting from any $\nu_0 \in 2\Z+1$.\par
%
%
\subsection{The Chaturvedi-Srinivasan Oscillator Algebra}
\subsubsection{Positive Values of the Deforming Parameter}
For the Chaturvedi-Srinivasan oscillator algebra ${\cal A}_1$~\cite{chaturvedi},
defined by Eqs.~(\ref{eq:2.1}) and~(\ref{eq:2.13}), we find from~(\ref{eq:3.8})
and~(\ref{eq:3.9}) that
\begin{equation}
  \lambda_n = \frac{q^{\nu_0+n}-1}{q-1} - c, \qquad \mbox{where } c\le
  \frac{q^{\nu_0}-1}{q-1},      \label{eq:4.16}
\end{equation}
is an increasing function of~$n$ for any positive~$q$.\par
%
%
Whenever $0<q<1$, $\lambda_n$ increases from~$-\infty$ to $(1-q)^{-1} - c$,
which is a positive constant as it follows from Eq.~(\ref{eq:4.16}). Hence, we only
get BFB~unirreps, characterized by
\begin{equation}
  \tnu_0 = \nu_0 + n_1, \qquad c = [\tnu_0]_q, \qquad \tlambda_n = q^{\tnu_0}
  [n]_q,       \label{eq:4.17}      
\end{equation}
where $n_1 \in \{0, -1, -2, \ldots\}$, showing that $\tnu_0$ may take any real
value.\par
%
%
Whenever $q>1$, $\lambda_n$ increases from $- (q-1)^{-1} - c$ to~$+\infty$, and
we have to distinguish between the cases where $- (q-1)^{-1} < c \le [\nu_0]_q$, and
$c\le - (q-1)^{-1}$. In the former, $\lambda_n$ may become negative, so that we
obtain BFB~unirreps similar to those defined in~(\ref{eq:4.17}). In the latter, on the
contrary, $\lambda_n$ is always positive, hence we get UB~unirreps.\par
%
%
In the limit where $q \to 1^-$ or~$1^+$, the UB~unirreps again diverge so that we
are only left with the BFB~ones, which go over into those of~${\cal A}_0$, as it was
the case for the Arik-Coon algebra.\par 
%
%
Comparing now the ${\cal A}_1$~unirreps with those of~$\tilde{\cal A}_0$, we
note that only the Fock-space representations of these algebras do coincide since
they are both characterized by $\tnu_0 = c = 0$, and $\tlambda_n = [n]_q$. The
remaining BFB~unirreps are however different and, more strikingly, the classically
singular representations appear for different $q$~values, namely $0<q<1$
for~$\tilde{\cal A}_0$ and $q>1$ for~${\cal A}_1$.\par
%
%
\begin{table}[htb]

\caption{Unirrep classification for the Chaturvedi-Srinivasan oscillator algebra.}

\vspace{1cm}
\begin{tabular}{lll}          
  \hline\\[-0.2cm] 
  $q$ & Type & Characterization \rule[-0.3cm]{0cm}{0.6cm}\\[0.2cm]
  \hline\\[-0.2cm] 
  $q>1$       & BFB & $\tnu_0\in\R$, $c=[\tnu_0]_q$, $\tlambda_n = q^{\tnu_0} [n]_q$
                       \\[0.2cm]
                 & UB   & $0\le\tnu_0<1$, $c\le - (q-1)^{-1}$, $\tlambda_n = [\tnu_0+n]_q
                       - c$\\[0.2cm]
  $0<q<1$   & BFB & $\tnu_0\in\R$, $c=[\tnu_0]_q$, $\tlambda_n = q^{\tnu_0}
                       [n]_q$\\[0.2cm]
  $-1<q<0$ & BFB & $\tnu_0\in2\Z$, $c=[\tnu_0]_q$, $\tlambda_n = q^{\tnu_0}
                       [n]_q$\\[0.2cm]
  $q=-1$    & FD  & $\tnu_0\in 2\Z$, $p=1$, $c=0$, $\tlambda_n = \left(1 - 
                       (-1)^n\right)/2$ \\[0.2cm]
                & UB   & $\tnu_0=0$, $c<0$, $\tlambda_n = - c + \left(1 - (-1)^n\right)/2$
                       \\[0.2cm]
  $q<-1$   & BFA  & $\tnu_0\in2\Z+1$, $c=[\tnu_0+1]_q$, $\tlambda_n = 
                       q^{\tnu_0+1} [n-1]_q$  \\[0.2cm]
  \hline 
\end{tabular}
\end{table}
%
%
\subsubsection{Negative Values of the Deforming Parameter}
For negative $q$~values and $\nu_0 \in \Z$, Eq.~(\ref{eq:4.16}) becomes
\begin{equation}
  \lambda_n = (-1)^{\nu_0+n+1} \frac{|q|^{\nu_0+n}}{1+|q|} + \frac{1}{1+|q|} - c, \qquad
  \mbox{where } c \le \frac{1-(-1)^{\nu_0}|q|^{\nu_0}}{1+|q|}.     \label{eq:4.18}
\end{equation}
For successive $n$~values, $\lambda_n$ oscillates around the constant $(1+|q|)^{-1}
- c$, which is positive for $\nu_0\in 2\Z$, but may be positive, null, or negative
for $\nu_0\in 2\Z+1$. To classify the unirreps, we have to distinguish between the
cases where $0<|q|<1$, $|q|>1$, and~$|q|=1$.\par
%
%
Whenever $0<|q|<1$, $|\lambda_n|$ decreases from~$+\infty$ to $|(1+|q|)^{-1} - c|$.
If $\nu_0 \in 2\Z$, $\lambda_n$ may become negative for some negative even
$n$~values. Unitarity then imposes that there exists some $n_1 \in \{0, -2, -4,
\ldots\}$ such that Eq.~(\ref{eq:4.9}) be satisfied. We therefore obtain
BFB~unirreps, characterized by
\begin{equation}
  \tnu_0 = \nu_0 + n_1 \in 2\Z, \qquad c = \frac{1-|q|^{\tnu_0}}{1+|q|}, \qquad
  \tlambda_n = |q|^{\tnu_0} \frac{1 - (-1)^n |q|^n}{1+|q|}.      \label{eq:4.19}
\end{equation}
If $\nu_0 \in 2\Z+1$, $c$ must satisfy the stronger condition $c\le (1+|q|)^{-1}$.
Then $\lambda_n$ may become negative for some negative odd $n$~values. Hence,
Eq.~(\ref{eq:4.9}) must be fulfilled for some $n_1 \in \{-1, -3, -5, \ldots\}$, so
that we get BFB~unirreps specified by~(\ref{eq:4.19}) again.\par
%
%
Whenever $|q|>1$, $|\lambda_n|$ increases from $|(1+|q|)^{-1} - c|$ to~$+\infty$.
Similar arguments show that Eq.~(\ref{eq:4.13}) must be fulfilled for some $n_2
\in \{2, 4, 6, \ldots\}$ or $n_2 \in \{1, 3, 5, \ldots\}$ according to whether $\nu_0
\in 2\Z$ or $\nu_0 \in 2\Z+1$. We therefore obtain BFA~unirreps, characterized by
\begin{equation}
  \tnu_0 = \nu_0 + n_2 - 1 \in 2\Z+1, \qquad c = \frac{1-|q|^{\tnu_0+1}}{1+|q|},
  \qquad \tlambda_n = |q|^{\tnu_0+1} \frac{1 + (-1)^n |q|^{n-1}}{1+|q|}.     
  \label{eq:4.20}
\end{equation}
\par
%
%
Finally, for $|q|=1$, hence for the oscillator algebra defined by~(\ref{eq:2.1}) and
$\left[a, \ap\right] = (-1)^N$, Eq.~(\ref{eq:4.18}) simply becomes
\begin{equation}
  \lambda_n = \frac{1}{2} \left(1 - (-1)^{\nu_0+n}\right) - c,     \label{eq:4.21}
\end{equation}
where $c\le 0$, or $c\le 1$, according to whether $\nu_0\in 2\Z$, or $\nu_0 \in
2\Z+1$. By reasoning as in Sec.~4.1.2, we find two types of unirreps, namely
FD~unirreps characterized by $c=0$, $p=1$, and any $\tnu_0 \in 2\Z$, and
UB~unirreps specified by $\tnu_0=0$, and any negative $c$~value.\par
%
%
\subsection{The Tamm-Dancoff Oscillator Algebra}
\subsubsection{Positive Values of the Deforming Parameter}
For the Tamm-Dancoff oscillator algebra $\tilde{\cal A}_1$~\cite{odaka}, defined
by Eqs.~(\ref{eq:2.1}) and~(\ref{eq:2.14}), we find from~(\ref{eq:3.8})
and~(\ref{eq:3.9}) that
\begin{equation}
  \lambda_n = q^{\nu_0+n-1} (\nu_0 + n - q c), \qquad \mbox{where } c\le q^{-1}
  \nu_0,              \label{eq:4.22}
\end{equation}
may become negative for $n < qc - \nu_0$, and any positive $q$~value. Hence, in the
present case, we only get BFB~unirreps, characterized by
\begin{equation}
  \tnu_0 = \nu_0 + n_1, \qquad c = q^{-1} \tnu_0, \qquad \tlambda_n = 
  q^{\tnu_0+n-1} n,       \label{eq:4.23}
\end{equation}
where $n_1 \in \{0, -1, -2, \ldots\}$, and therefore $\tnu_0 \in \R$.\par
%
%
Such unirreps were already found before~\cite{odaka}. The fact that $\lim_{n\to
+\infty} \tlambda_n = 0$, for $0<q<1$, explains the name given to the algebra, and
referring to the idea of a high-energy cutoff proposed in the context of field
theory~\cite{tamm}.\par
%
%
\begin{table}[htb]

\caption{Unirrep classification for the Tamm-Dancoff oscillator algebra.}

\vspace{1cm}
\begin{tabular}{lll}          
  \hline\\[-0.2cm] 
  $q$ & Type & Characterization \rule[-0.3cm]{0cm}{0.6cm}\\[0.2cm]
  \hline\\[-0.2cm] 
  $0<q\ne1$       & BFB & $\tnu_0\in\R$, $c=q^{-1}\tnu_0$, $\tlambda_n = 
             q^{\tnu_0+n-1} n$ \\[0.2cm]
  \hline 
\end{tabular}
\end{table}
%
%
\subsubsection{Negative Values of the Deforming Parameter}
For negative $q$~values and $\nu_0 \in \Z$, Eq.~(\ref{eq:4.22}) becomes
\begin{equation}
  \lambda_n = (-1)^{\nu_0+n+1} |q|^{\nu_0+n+1} (\nu_0 + n + |q| c),    \label{eq:4.24}
\end{equation}
where $c \le - \nu_0 |q|^{-1}$ or $c \ge - \nu_0 |q|^{-1}$ according to whether
$\nu_0 \in 2\Z$ or $\nu_0 \in 2\Z+1$. Since $\lambda_n$ can vanish for at most
one integer $n$~value, and it oscillates around zero in the intervals $(-\infty,
-\nu_0 - |q| c)$ and $(-\nu_0 - |q| c, +\infty)$, it is obvious that the conditions of
Proposition~1 cannot be fulfilled so that no unirrep can exist for negative $q$
values.\par
%
%
We have therefore established that contrary to the remaining deformed oscillator
algebras considered in the present paper, the Tamm-Dancoff oscillator algebra has a
single class of unirreps.\par
%
%
\section{Conclusion}
In the present paper, we developed the representation theory of deformed oscillator
algebras, defined in terms of an arbitrary function of the number operator~$N$. We
showed that the classification of their unirreps can be most easily performed in
terms of the eigenvalues of a Casimir operator~$C$. Under the assumption that the
spectrum of~$N$ is discrete and nondegenerate, we proved that the unirreps may
fall into one out of four classes (BFB, BFA, FD, UB) according to the nature of that
spectrum, and that bosonic, and fermionic or parafermionic Fock-space
representations may occur as special cases of BFB and FD~unirreps,
respectively.\par
%
%
We did also carry out the unirrep classification in detail for some deformed
oscillator algebras, which can be derived from the boson one by the recursive
minimal deformation procedure of Katriel and Quesne~\cite{katriel}, namely the
Arik-Coon~\cite{arik,kuryshkin}, Chaturvedi-Srinivasan~\cite{chaturvedi}, and
Tamm-Dancoff~\cite{odaka} oscillator algebras. For all of them, we considered
both positive and negative values of the deforming parameter, which constitutes a
distinctive feature of the present study as compared with some previous
ones~\cite{kulish,rideau,chaichian,kosinski}.\par
%
%
We showed that all the known unirreps, in particular the bosonic Fock-space
representations, can be recovered in our classification scheme, and that in addition,
many new unirreps make their appearance. We actually provided some examples for
each of the four unirrep classes, although in the FD~case, only two-dimensional
unirreps were encountered. Higher-dimensional FD~unirreps do however arise
for some known deformed oscillator algebras~\cite{cqb}.\par
%
%
We also illustrated both the effects of minimal deformation and of the
quommutator-commutator transformation of the recursive procedure on
non-Fock-space representations.\par
%
%
Applications of deformed oscillator algebras have been restricted up to now to
their Fock-space representations. Whether non-Fock-space representations, such as
those constructed in this paper, may have some useful applications remains an
interesting open question.\par
%
%

\end{document}